# On Noneven Digraphs and Symplectic Pairs

Chjan C. Lim

David A. Schmidt

Department of Mathematics

Rensselaer Polytechnic Institute

# 1 Abstract

A digraph $D$ is called **noneven** if it is possible to assign weights of 0,1 to its arcs so that $D$ contains no cycle of even weight. A noneven digraph $D$ corresponds to one or more nonsingular sign patterns.

Given an $n \times n$ sign pattern $H$, a **symplectic pair** in $Q(H)$ is a pair of matrices $(A, D)$ such that $A \in Q(H)$, $D \in Q(H)$, and $A^T D = I$. An unweighted digraph $D$ allows a matrix property $P$ if at least one of the sign patterns whose digraph is $D$ allows $P$.

In [1], Thomassen characterized the noneven, 2-connected symmetric digraphs (i.e., digraphs for which the existence of arc $(u, v)$ implies the existence of arc $(v, u)$). In the first part of our paper, we recall this characterization and use it to determine which strong symmetric digraphs allow symplectic pairs.

A digraph $D$ is called **semi-complete** if, for each pair of distinct vertices $(u, v)$, at least one of the arcs $(u, v)$ and $(v, u)$ exists in $D$. Thomassen, again in [1], presented a characterization of two classes of strong, noneven digraphs: the semi-complete class and the class of digraphs for which each vertex has total degree which exceeds or equals the size of the digraph. In the second part of our paper, we fill a gap in these two characterizations and present and prove correct versions of the main theorems involved. We then proceed to determine which digraphs from these classes allow symplectic pairs.

# 2 Introduction

A **sign pattern** is an $m \times n$ matrix $H$ whose entries are 0,1 or -1. We say that an $m \times n$ matrix $A \in Q(H)$ (the sign-class of $H$) if $sgn(A_{ij}) = sgn(H_{ij})$ for all pairs of indices $(i, j)$,

where
$$sgn(x) = \begin{cases} -1, & x < 0; \\ 0, & x = 0; \\ 1, & x > 0. \end{cases}$$

An $n \times n$ sign pattern $H$ is called **sign-nonsingular** (or sometimes an **L-Matrix**) if for all matrices $A \in Q(H)$, $det(A) \neq 0$. $H$ is **combinatorially singular** if $det(A) = 0$ for all $A \in Q(H)$. Sign-nonsingular matrices have been studied extensively in the context of the sign-solvability problem [2],[3].

The **weighted digraph** of an $n \times n$ sign pattern $H$ is the directed graph $D(H)$ with $n$ vertices and the arc set $\{(i,j) : i \neq j, H_{ij} \neq 0\}$, with each arc assigned the weight 1 if $H_{ij} < 0$ and the weight 0 otherwise. Since the diagonal entries of $H$ are not reflected in $D(H)$, it is convenient when possible to restrict consideration to **negative-diagonal** sign patterns, i.e. sign patterns $H$ for which $H_{ii} = -1$ for all indices $i$.

Our digraph terminology follows that of [1]. A digraph $D$ is **strongly connected** (or **strong**) if, for any ordered pair of vertices $(u, v)$ from $D$, there is a directed path beginning at $u$ and terminating at $v$. $D$ is **strongly $k$-connected** if the removal of $k - 1$ or fewer vertices from $D$ results in a strong digraph. Vertex $u$ **dominates** vertex $v$ in $D$ if the arc $(u, v)$ exists in $D$. The **indegree** of vertex $v$ ($\text{id}(v)$) is the number of vertices that dominate $v$, and the **outdegree** of $v$ ($\text{od}(v)$) is the number of vertices dominated by $v$. The **degree** of $v$ ($\text{d}(v)$) is the sum of $v$'s indegree and outdegree. Let $D$ be a directed graph of size $n$ with vertex set $V$. A finite sequence $\{d_1, d_2, ..., d_n\}$ is an **indegree sequence** for $D$ if there is an ordering $\{v_i\}$, $i = 1, 2, ..., n$, of $V$ such that $\text{id}(v_i) = d_i$ for $i = 1, 2, ..., n$. An **outdegree sequence** for $D$ is defined analogously. $D$ is called **semi-complete** if for any pair of vertices $(u, v)$ in $D$, either $u$ dominates $v$ or vice versa. The undirected cycle of length $k > 2$ is denoted $C_k$, with $C_2$ denoting a single edge joining 2 vertices.

A weighted digraph $D$ is called **even** if some directed cycle in $D$ has even total weight,

where the weight of a directed path or cycle in $D$ is obtained by summing the weights of the arcs contained in that path or cycle. If $D$ is not even, it is called **noneven**. A negative-diagonal sign-nonsingular matrix has the property [4] that its weighted digraph is noneven.

An unweighted digraph $D$ is called **even** if $D$ contains an even cycle whenever the arcs of $D$ are assigned weights of 0 or 1. If there is at least one weight assignment for the arcs of $D$ that results in a noneven weighted digraph, $D$ is called **noneven**. In view of the above, it follows that any noneven digraph $D$ corresponds to at least one, and perhaps several, negative-diagonal sign-nonsingular patterns.

A digraph $D$ is called **symmetric** if the existence of arc $(u, v)$ in $D$ implies the existence of the arc $(v, u)$. A negative-diagonal sign pattern whose digraph is symmetric is called **combinatorially symmetric**. Every undirected graph $G$ gives rise to a corresponding symmetric digraph which we shall call $G^*$. An undirected graph $G$ will be called noneven if $G^*$ is noneven.

A **subdivision** of an arc $(u, v)$ in a digraph $D$ is obtained by removing $(u, v)$, adding a new vertex $w$, and adding the arcs $(u, w)$ and $(w, v)$. To **split** a vertex $v$ in $D$ means to replace $v$ with two vertices $v_1$ and $v_2$ and an arc $(v_1, v_2)$, such that all arcs dominating (respectively, dominated by) $v$ in $D$ dominate $v_1$ (respectively, are dominated by $v_2$) in the resulting digraph. The result of splitting vertices and subdividing arcs in the double cycle $C_k^*$ is called a **weak $k$-double-cycle**. Seymour and Thomassen showed [5] that a digraph $D$ is even if and only if $D$ contains a weak $k$-double-cycle for some odd $k$.

Using terminology from [6], let $P$ be a property that a real, square matrix may or may not possess. The sign pattern $H$ is said to **allow** property $P$ if there exists a matrix $A \in Q(H)$ such that $A$ has property $P$. $H$ is said to **require** $P$ if all matrices $A$ in $Q(H)$ have property $P$. Two sign patterns $G$ and $H$ are **sign-equivalent** (written $G \sim H$) if $G$ can be obtained from $H$ by negating some of its rows and columns and then permuting the rows and columns.

Given two negative-diagonal sign patterns $G$ and $H$, we will call their weighted digraphs $D(G)$ and $D(H)$ sign-equivalent if $G \sim H$. A matrix property $P$ is called a **sign-** **equivalence class property** if a sign pattern $H$ allowing $P$ implies that $P$ is allowed by all sign patterns that are sign-equivalent to $H$. The problem of characterizing the patterns allowing such a property is simplified because the analysis can be carried out using only negative-diagonal sign patterns, since any pattern that is not combinatorially singular is sign-equivalent to some negative-diagonal pattern.

A weighted digraph $D$ is said to allow matrix property $P$ if the negative-diagonal sign pattern corresponding to $D$ allows $P$. An unweighted digraph $D$ is said to allow $P$ if some negative-diagonal sign pattern whose digraph is $D$ allows $P$.

A result of Brualdi and Shader [7] states that any two negative-diagonal nonsingular sign patterns with the same unweighted noneven digraph $D$ are sign-equivalent. This result implies that two weighted noneven digraphs $D_1$ and $D_2$ are sign-equivalent if and only their indegree and outdegree sequences coincide. It therefore allows us to classify completely the negative-diagonal sign-nonsingular patterns allowing a sign-equivalence class property $P$ using only unweighted digraphs.

Given an $n \times n$ sign pattern $H$, a **symplectic pair** in $Q(H)$ is a pair of matrices $(A, D)$ such that $A \in Q(H), D \in Q(H)$, and $A^T D = I$. Symplectic pairs are a pattern-generalization of orthogonal matrices which arise from a special symplectic matrix found in $n$-body problems in celestial mechanics, and have been discussed previously in [8],[9]. In particular, it was shown in [9] by the authors that the symplectic pair property is a sign-equivalence class property. In view of this fact, the characterizations of sign-nonsingular patterns allowing symplectic pairs discussed in this paper are carried out in terms of unweighted digraphs. It should be understood that if a noneven digraph $D$ is shown to allow symplectic pairs, it follows that any negative-diagonal sign-nonsingular pattern whose digraph is $D$ (and any sign pattern that is sign-equivalent to such a pattern) also allows symplectic pairs.

# 3 Overlap Numbers, Digraphs and Symplectic Pairs

We first discuss a necessary condition for a sign pattern to allow symplectic pairs. Let $H$ be an $n \times n$ sign pattern. Given disctinct indices $i$ and $j$, define the **signed row overlap numbers** $N(H)_{ij}^+$ (respectively, $N(H)_{ij}^-$) as the number of column indices $k$ for which $H_{ik}H_{jk} > 0$ (respectively, $< 0$). The following necessary condition for a sign pattern $H$ to allow symplectic pairs follows directly from the requirement that the sign pattern $H$ permit orthogonality of row vectors for $A \in Q(H)$:

**Proposition 3.1:** Let $H$ be an $n \times n$ sign pattern. If $N(H)_{ij}^+ N(H)_{ij}^- = 0$ for some pair of distinct indices and $H$ allows symplectic pairs, then $N(H)_{ij}^+ = N(H)_{ij}^- = 0$.

We will need a digraphical version of this necessary condition. Assuming for the moment that a vertex in the digraph $D(H)$ dominates itself if $H_{ii} \neq 0$, Proposition 3.1 can be restated as follows:

**Proposition 3.1** (digraphical version): Let $D(H)$ be the weighted digraph of the sign pattern $H$. Given a pair of distinct vertices $v$ and $w$ from $D(H)$, let $U_{vw}$ be the set of all vertices from $D(H)$ which are dominated by both $v$ and $w$. Then $U_{vw}$ is the disjoint union of the two sets $U_{vw}^+$ and $U_{vw}^-$, the first consisting of the vertices $u$ for which the sum of the weights of the arcs $(v, u)$ and $(w, u)$ is even, and the second consisting of those for which this sum is odd. If there exists a pair of distinct vertices $v$ and $w$ for which exactly one of the sets $U_{vw}^+$ and $U_{vw}^-$ is empty, $D(H)$ does not allow symplectic pairs.

Figure 1 shows a negative-diagonal sign pattern $H$ (with its weighted digraph $D(H)$) which fails to allow symplectic pairs by Proposition 3.1. Note that for a negative-diagonal pattern $H$, the self-arcs of $D(H)$ inferred here all have unit weight.

$$H = \begin{pmatrix} -1 & 1 & 0 & -1 \\ 0 & -1 & 0 & 1 \\ 0 & -1 & -1 & 1 \\ 0 & -1 & 0 & -1 \end{pmatrix}$$

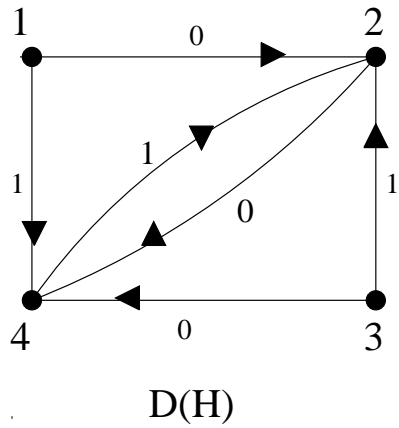

Figure 1: $H$ and $D(H)$, Where $H$ Fails to Allow Symplectic Pairs

# 4 Symmetric Noneven Digraphs and Symplectic Pairs

A digraph $D$ is called **symmetric** if the existence of arc $(u, v)$ in $D$ implies the existence of the arc $(v, u)$. A negative-diagonal matrix whose digraph is symmetric is called **combinatorially symmetric**. Every undirected graph $G$ gives rise to a corresponding symmetric digraph which we shall call $G^*$. An undirected graph $G$ will be called noneven if $G^*$ is noneven.

In [1], Thomassen characterized the 2-connected noneven symmetric digraphs. Construct the graph $G$ as follows. Starting with $C_4$, choose any edge $(u, v)$ and add the vertices $u'$ and $v'$, and add the edges $(u, u')$, $(u', v')$ and $(v', v)$. Repeat this procedure as desired. The resulting 2-connected graph $G$ is called a $C_4$-cockade. Thomassen showed that a 2-connected graph $G$ is noneven if and only if $G$ is a subgraph of a $C_4$-cockade, and also proved a useful equivalent condition, which was proven independently by Harary, et al [10]:

**Proposition 4.1:** The 2-connected graph $G$ is even if and only if $G = C_n$, $n$ odd, or $G$ contains two vertices joined by three internally disjoint paths, one of which has even length.

We shall call a cycle $C$ in a graph $G$ **reduced** if the subgraph of $G$ induced by the vertices from $C$ is $C$ itself. The following facts are direct consequences of Proposition 4.1. If $G$ is a 2-connected subgraph of a $C_4$-cockade, $G$ must contain at least one reduced cycle of length 4 or greater. In particular, a 2-connected subgraph $G$ of a $C_4$-cockade is itself a $C_4$-cockade if and only if all of its reduced cycles have length 4.

We now consider the problem of determining which of the strong, symmetric noneven digraphs allow symplectic pairs.

**Theorem 4.2:** Let $G^*$ be a strong, symmetric noneven digraph with underlying undirected graph $G$. If $G^*$ allows symplectic pairs, then either $G^* = C_2^*$ or $G^* = C_4^*$.

**Proof:** If $G^* = C_2^*$, then the sign-nonsingular pattern

$$H = \begin{pmatrix} -1 & 1 \\ -1 & -1 \end{pmatrix},$$

whose unweighted digraph is $G^*$, allows, indeed requires symplectic pairs (see Chapter Three). If $G^* = C_4^*$, then the negative-diagonal sign-nonsingular pattern

$$H = \begin{pmatrix} -1 & -1 & 0 & -1 \\ 1 & -1 & -1 & 0 \\ 0 & 1 & -1 & -1 \\ 1 & 0 & 1 & -1 \end{pmatrix},$$

whose digraph if $G^*$, allows symplectic pairs and also allows the more restrictive property of orthogonality (see [11]).

Now suppose $G^*$ is neither $C_2^*$ nor $C_4^*$. If $G$ has connectivity one, then $G$ has a cut vertex $v$, and there exist vertices $u$ and $w$ (distinct from $v$) in $G$ such that every $uw$-path contains $v$. Consider an arbitrary $uw$-path $P$ (at least one of which must exist since $G$ is connected), and call the unique pair of vertices that are adjacent to $v$ in $P$ $u'$ and $w'$. In $G^*$, both $u'$ and $w'$ dominate $v$, but cannot both dominate any other vertex of $G$ since this would result in a $uw$ path in $G$ not containing $v$. Thus $G$ fails to allow symplectic pairs by Proposition 3.1.

If $G$ is 2-connected, $G$ must be a subgraph of a $C_4$-cockade. Suppose first that $G$ is a $C_4$-cockade. Since $G \neq C_4$, $G$ must have the form shown in Figure 2, where the subgraph $\overline{G}$ of $G$ is also a $C_4$-cockade.

Note that $u'$ and $u''$ both dominate $u$ in $G^*$, and that the only other vertices dominated by $u'$ in $G^*$ are itself and $v'$. But $u''$ dominates neither $u'$ nor $v'$. Thus $G$ fails to allow symplectic pairs by Proposition 3.1.

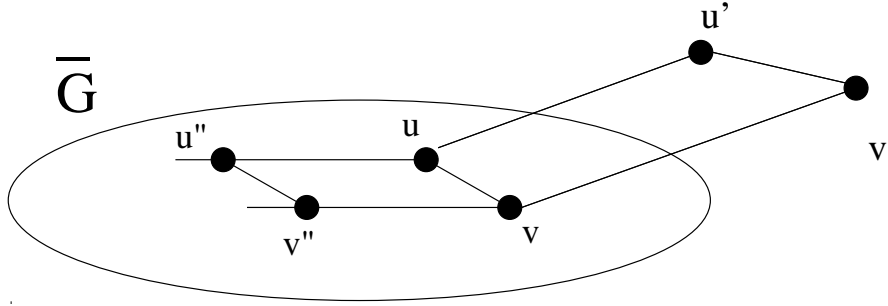

Figure 2: A $C_4$-Cockade Properly Containing $C_4$

Now suppose that $G$ is not a $C_4$-cockade, which means that $G$ contains a reduced cycle $C$ of length $p \geq 6$. Consider three vertices $u, v$, and $w$ from $C$ (as shown in Figure 3) such that $u$ and $v$ both dominate $w$ in $G^*$. Since $C$ is a reduced cycle, $u$ and $v$ cannot both be adjacent to any other vertex in $C$, since otherwise $C$ would have length four. Suppose $u$ and $v$ are both adjacent to a vertex $x$ not in $C$, as shown in Figure 3. Then there exist three internally disjoint paths joining $u$ and $v$, namely the two paths from $u$ to $v$ on $C$ and the path $uxv$. Furthermore, this last path has length two. Thus by Proposition 4.1, $G$ is even, which is a contradiction. Thus no such vertex $x$ exists in $G$, and $G$ fails to allow symplectic pairs by Proposition 3.1. ∎

Note that a symmetric digraph $G^* \neq C_2^*$ such that $G$ has connectivity one fails to allow symplectic pairs whether $G^*$ is noneven or not.

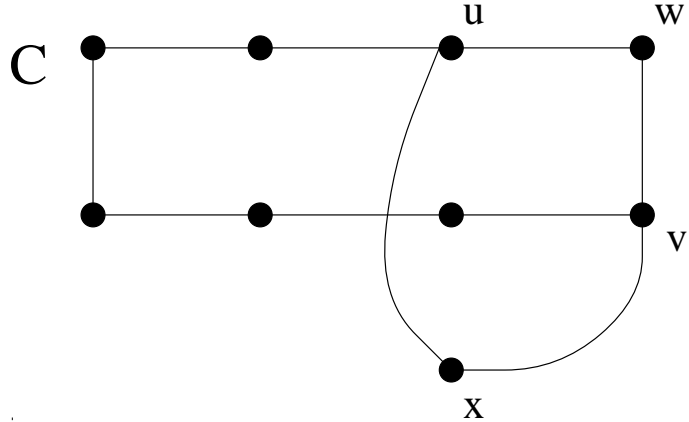

Figure 3: A Reduced Cycle of Length $p \geq 6$

## 5 Dense Noneven Digraphs

We shall call a digraph $D$ **dense** if $D$ has one of the properties studied by Thomassen with regard to the even cycle problem in [1]; i.e., either $D$ is semi-complete, or $D$ is such that the total degree of each vertex exceeds or equals the size of $D$. In the study of dense noneven digraphs, an important role is played by the **extended caterpillars**. An extended caterpillar is a semi-complete digraph constructed as follows. Consider a simple undirected path $P = x_1 x_2 ... x_k$. With each vertex $x_j, 2 \leq j \leq k-1$, associate a (possibly empty) set of vertices $\{x_{j_i}\}$, all of which are adjacent to $x_j$. Call the resulting tree $T$. Call the path $P$ the **backbone** of $T$, and call a particular backbone vertex $x_j$ together with its associated vertices $\{x_{j_i}\}$ (if any) a **blossom** of $T$. Note that the ordering of the backbone vertices of $T$ induces (through blossom membership) an ordering on all the vertices of $T$; in particular, we will write $u > v$ if $u$ is in the blossom associated with backbone vertex $x_i$, $v$ is in the blossom associated with backbone vertex $x_j$, and $i > j$. We now form the extended caterpillar $E$ by adding to $T^*$ all arcs $(u,v)$ not present in $T^*$ such that $u > v$, and exactly one arc joining any two non-backbone vertices in a given blossom, so that each blossom contains a transitive tournament on its vertex set. It is clear that any extended caterpillar is dense in both of the senses discussed above, and it is easy to show that an extended caterpillar is noneven. An example of an extended caterpillar of size six is shown in Figure 4.

Extended caterpillars and their strong subdigraphs have connectivity one. Also, a strong subdigraph of an extended caterpillar $E$ must contain a backbone path of the form $x_j...x_m$, where $1 \leq j \leq m \leq k$ and $\mathrm{id}(x_m) = \mathrm{od}(x_j) = 1$. In [1], Thomassen gave the following

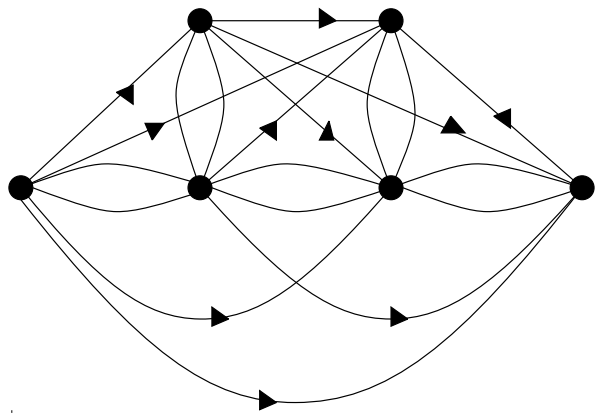

Figure 4: An Extended Caterpillar of Size Six

characterization of the strong semi-complete noneven digraphs:

> Let $D$ be a strong semi-complete digraph. Then $D$ is noneven if and only if $D$ is a subdigraph of an extended caterpillar.

This characterization, while valid if $D$ has connectivity one, neglects the 2-connected digraph (which we shall call $W_4$ due to its relationship with a certain wheel graph [11]) shown in Fig-ure 5. $W_4$ is strong, semi-complete and noneven, yet cannot be a subdigraph of an extended caterpillar due to its 2-connectedness. We now state and prove a correct characterization for the 2-connected case. In this proof, as well as in the Proof of Theorem 5.2 that follows, we have used ideas from the corresponding proofs in [1].

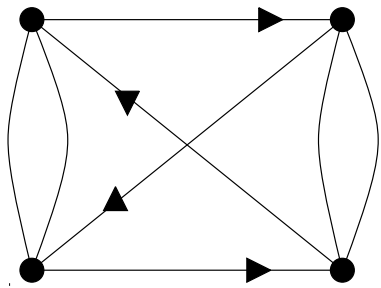

Figure 5: An Extended Caterpillar of Size Six

**Theorem 5.1:** Let $D$ be a 2-connected semi-complete noneven digraph. Then $D = W_4$.

**Proof:** If $n = 4$, then $D = W_4$, since a 2-connected noneven digraph of size four must have exactly eight arcs with each vertex of indegree and outdegree two, and $W_4$ is the only such semi-complete digraph.

The proof proceeds by induction. Suppose $D$ has size $n \geq 5$. We first consider the case where $D$ contains a vertex $x$ for which $D - x$ is 2-connected. By the induction hypothesis, $n = 5$ and $D - x = W_4$ in this case. Since $D$ is 2-connected and semi-complete, $D$ must contain as a subdigraph one of the digraphs shown in Figure 6. Both of these digraphs contain a weak 3-double-cycle, and thus in either case, $D$ is even, which is a contradiction.

Now suppose that, for every vertex $x$ in $D$, $D - x$ has connectivity one, and is thus a subdigraph of an extended caterpillar. Let $\{u_i\}$ be the set of vertices from $D - x$ with outdegree one, and $\{v_i\}$ the set of vertices from $D - x$ with indegree one. (Since any subdigraph of an extended caterpillar has at least one vertex of each of these types, both of these sets are nonempty, and since $D - x$ has size four or larger, semi-completeness implies that these two sets do not intersect.) Since $D$ is 2-connected, $x$ dominates each vertex in $\{v_i\}$ and is dominated by each vertex in $\{u_i\}$. Given a vertex $u \in \{u_i\}$, construct the digraph $D'_u$ by adding to $D - x$ all arcs from $u$ to vertices in $\{v_i\}$. Given a vertex $v \in \{v_i\}$, construct the digraph $D'_v$ by adding to $D - x$ all arcs from vertices in $\{u_i\}$ to $v$. As noted in [1], such a $D'_u$ or $D'_v$ is still noneven, since the arc $(u, v_i)$ (respectively, $(u_i, v)$) can be weighted identically to the path $uxv_i$ (respectively, $u_ixv$).

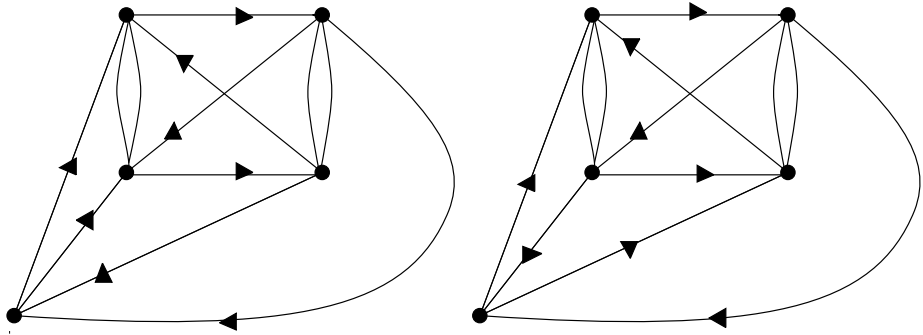

Figure 6: See Proof of Theorem 5.1

Suppose that for all $u \in \{u_i\}$ and all $v \in \{v_i\}$, $D'_u$ and $D'_v$ have connectivity one, and fix $u \in \{u_i\}$, so that $D'_u$ is a subdigraph of an extended caterpillar. Since $u$ has outdegree one in $D - x$, $u$ dominates at most one vertex from $\{v_i\}$. If $u$ dominates no vertex from $\{v_i\}$, $D'_u$ has no vertex of indegree one, a contradiction. Thus $u$ dominates exactly one vertex from $\{v_i\}$. Given any two vertices from $\{u_i\}$, one of them must dominate the other since $D - x$ is semi-complete. But this is a contradiction since each vertex of outdegree one must dominate a vertex of indegree one. This shows that $D - x$ contains exactly one vertex $u$ of outdegree one. A similar argument shows that $D - x$ contains exactly one vertex $v$ of indegree one. Since $u$ dominates $v$ and $D - x$ is a subdigraph of an extended caterpillar, we conclude that $D - x$ contains only two vertices, which is a contradiction.

Now suppose there is a $u \in \{u_i\}$ such that $D'_u$ is 2-connected, which implies that $n = 5$ and $D'_u = W_4$. Then $u$ cannot dominate any $v \in \{v_i\}$ in $D - x$, otherwise $D'_u$ would contain a vertex of indegree one, contradicting its 2-connectedness. Thus, by semi-completeness, each $v \in \{v_i\}$ dominates $u$. If there were more than one vertex of indegree one in $D - x$, then $D'_u$ would contain two 2-cycles incident at the same vertex, which $W_4$ does not. Thus there is a unique $v \in \{v_i\}$, and so $D - x$ is obtained from $W_4$ by removing an arc from one of the double cycles. A similar argument yields the same form for $D - x$ in the case where there is a $v \in \{v_i\}$ such that $D'_v$ is 2-connected. Now since $x$ dominates $v$ and is dominated by $u$, $D$ must be isomorphic to one of the digraphs appearing in Figure 7. Since both contain a weak 3-double-cycle, we again obtain the contradiction that $D$ is even, completing the proof. ∎

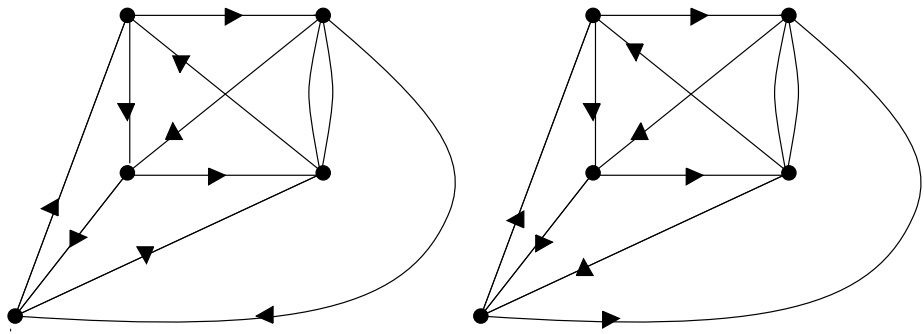

Figure 7: See Proof of Theorem 5.1

The second type of dense digraph is one in which the total degree of each vertex equals or exceeds the size of $D$. The characterization of the strong, noneven digraphs in this class from [1] is as follows:

> If $D$ is a strong, noneven digraph with $n$ vertices and of minimum degree at least $n$, then $D$ is a subdigraph of an extended caterpillar or else $D$ is isomorphic to $C_4^*$.

This result is correct if $D$ has connectivity one, but neglects $W_4$, a 2-connected noneven digraph which satisfies $\mathrm{d}(v) = n$ for all its vertices. A correct version for the 2-connected case follows:

**Theorem 5.2:** Suppose $D$ is a 2-connected, noneven digraph of size $n$ such that $\mathrm{d}(v) \geq n$ for all its vertices $v$. Then either $D = W_4$ or $D = C_4^*$.

**Proof:** If $n = 4$, either $D = W_4$ or $D = C_4^*$, as these are the only two noneven 2-connected digraphs of size four.

Proceeding by induction, suppose $D$ has size $n \geq 5$. By a result of Thomassen and Haagkvist [12], $D$ must contain a cycle $x_1, x_2, ..., x_{n-1}$ which misses exactly one of its vertices. Call the vertex not in this cycle $x$. Since $D$ is noneven, it follows that there is at most one 2-cycle incident at $x$, say from $x_{n-1}$. Since $x$ has degree $n$, there must be exactly one 2-cycle incident at $x$, and every vertex of $D - x$ either dominates or is dominated by $x$. If there exist $i, j$ such that $1 \leq i < j \leq n - 2$ and the arcs $(x, x_i)$ and $(x_j, x)$ are present, then $D$ contains a subdivision of the digraph shown in Figure 8, and is thus even, which is a contradiction. Thus there is a $k$, $1 \leq k \leq n - 2$ such that the vertices $x_1, x_2, ..., x_k$ dominate $x$, and the vertices $x_{k+1}, x_{k+2}, ..., x_{n-2}$ are dominated by $x$.

Construct the digraph $D'$ from $D - x$ as follows. If each vertex of $D - x$ has degree at least $n - 1$, then let $D' = D - x$. If $D - x$ has a vertex of degree $n - 2$, there must either be a vertex $x_i$, $1 \leq i \leq k$, such that $x_i$ does not dominate $x_{n-1}$, or a vertex $x_j$, $k + 1 \leq j \leq n - 2$, such that $x_{n-1}$ does not dominate $x_j$. Form $D'$ by adding to $D - x$ either the arc $(x_i, x_{n-1})$ or the arc $(x_{n-1}, x_j)$. In either case, $D'$ is noneven, since the arc $(x_i, x_{n-1})$ (respectively, $(x_{n-1}, x_j)$) can be weighted identically to the path $x_i x x_{n-1}$ (respectively $x_{n-1} x x_j$) in $D$. Furthermore, each vertex in $D'$ has degree at least $n-1$.

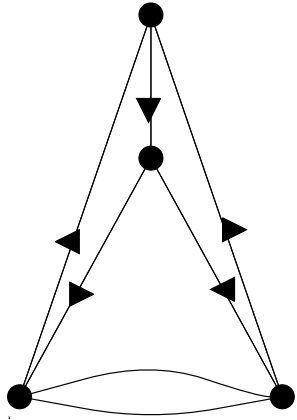

Figure 8: See Proof of Theorem 5.2

We first consider the case where $D'$ is 2-connected. By the induction hypothesis, $n = 5$ and $D'$ is either $W_4$ or $C_4^*$. If $D' = D - x = C_4^*$, then $D$ must contain the digraph shown in Figure 9. No matter how the undirected edges are oriented (so long as $D$'s 2-connectedness is required) it is easily confirmed that $D$ must contain a weak 3-double-cycle, contradicting the fact that $D$ is noneven. If $D' = D - x = W_4$, then $D$ must contain one of the digraphs shown in Figure 6, and $D$ is again even. If $D'$ was obtained from $D - x$ by the addition of an arc, then $D - x$ must have had two vertices of degree three, which is a contradiction. Thus $D'$ has connectivity one, and is therefore a subdigraph of an extended caterpillar. Let the basic path of $D'$ be $z_1 z_2 ... z_m$, so that od $(z_1) =$ id $(z_m) = 1$.

**Claim:** The only vertex of outdegree one in $D'$ is $z_1$, and the only vertex of indegree one in

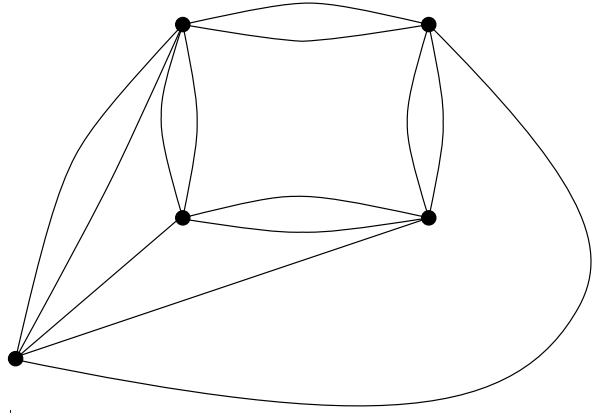

Figure 9: See Proof of Theorem 5.2

$D'$ is $z_m$.

**Proof of Claim:** $D'$ contains the hamiltonian cycle $x_1 x_2 ... x_{n-1} x_1$. Since $\text{od}(z_1) = 1$ and $\text{d}(z_1) \geq n - 1$, every vertex in $D'$ other than $z_1$ dominates $z_1$. Thus the only vertex in $D'$ other than $z_1$ that could possibly have outdegree one is the vertex $z_0$ which precedes $z_1$ in the hamiltonian cycle. But if $\text{od}(z_0) = 1$ then $z_1$ dominates $z_0$ so that $z_0 = z_2$. This is impossible since $n \geq 5$. Thus $z_1$ is the only vertex of outdegree one in $D'$. A similar argument shows that $z_m$ is the only vertex of indegree one in $D'$.

**Claim:** Suppose $z_m \neq x_{n-1}$. Then the only vertices of indegree one in $D - x$ are $z_m$ and possibly $x_{n-1}$.

**Proof of Claim:** If $D' = D - x$ or $D'$ is obtained from $D - x$ by the addition of the arc $(x_j, x_{n-1}), 1 \leq j \leq k$, then there is nothing to prove. So suppose $D'$ is obtained from $D - x$ by adding the arc $(x_{n-1}, x_j), k + 1 \leq j \leq n - 1$. Note that $x_j \neq z_m$, since otherwise $\text{id}(z_m) = 0$ in $D - x$, a contradiction since $D - x$ is strong. If $\text{id}(x_j) = 1$ in $D - x$, then $x_j$ must dominate every other vertex in $D - x$, since $\text{d}(x_j) \geq n - 1$. Thus $x_j$ dominates $z_m$, so that $x_j = z_{m-1}$. Since $z_m \neq x_{n-1}$, $\text{d}(z_m) \geq n - 1$, and $z_m$ dominates $z_{m-1}$. Noting that $z_{m-2}$ dominates $z_{m-1}$ as well, we see that $x_j$ cannot have indegree one.

We can assume without loss of generality that $z_1 \neq x_{n-1}$, since if this were the case, we could consider the digraph obtained from $D - x$ by reversing every arc, a digraph which is noneven if and only if $D - x$ is.

**Claim:** In $D - x$, $z_1$ does not dominate $x_{n-1}$.

**Proof:** If $z_1$ dominates $x_{n-1}$, then $x_{n-1} = z_2$. Furthermore, since the entire backbone of $D - x$ must lie on the hamiltonian cycle $x_1 x_2 ... x_{n-1} x_1$, $z_1 = x_{n-2}$. Thus, since $z_1$ has outdegree one, $z_1$ dominates $x$ in $D$, and so all vertices $x_i, 1 \leq i \leq n - 3$ dominate $x$ in $D$ as well. Thus, in $D$, $\text{od}(x) = 1$, a contradiction since $D$ is 2-connected.

Now form the digraph $D''$ from $D - x$ by adding all arcs of the form $(z_1, y)$ where $y$ is dominated by $x$ in $D$. Since $z_1$ does not dominate $x_{n-1}$ in $D - x$, and $x$ dominates $x_{n-1}$ in $D$, each vertex $v$ of $D''$ satisfies $\text{d}(v) \geq n - 1$. Also, $D''$ is noneven since the arc $(z_1, y)$ in $D''$ can be weighted the same as the path $z_1 x y$ in $D$. Thus $D''$ is a subdigraph of an extended caterpillar. Now, in $D - x$, $z_1$ doesn't dominate $z_m$, since $n \geq 5$. Furthermore, $x$ dominates both $x_{n-1}$ and $x_m$ in $D$. Thus, if $z_m \neq x_{n-1}$, $D''$ has no vertex of indegree one, which is a contradiction. If $z_m = x_{n-1}$, then either $z_m$ is the only vertex of indegree one in $D - x$ (in which case the preceding argument leads to a contradiction), or $z_1$ is the only vertex of outdegree one in $D - x$. Since $D''$ is obtained from $D - x$ by adding at least one arc incident from $z_1$, the latter case implies that $D''$ has no vertex of outdegree one, which is a contradiction that completes the proof. ∎

# 6 Dense Noneven Digraphs and Symplectic Pairs

In order to simplify our characterization of the dense noneven digraphs allowing symplectic pairs, we first explore a second necessary condition for a strong digraph to allow symplectic pairs.

Note that each non-zero term in the expansion of $K_A[i,j]$ corresponds to a signed path from $j$ to $i$ in $D(A)$ (see Appendix A). Recalling the definition of maximality of a sign-nonsingular pattern, we shall refer to a noneven digraph as maximal if the addition of any arc results in an even digraph.

If a sign pattern $H$ is sign-nonsingular, we can use Lemma 2.3-1 to construct the following necessary condition for $H$ to allow symplectic pairs:

**Proposition 6.1:** Let the sign pattern $H$ be sign-nonsingular. Then $H$ can allow symplectic pairs only if the following two conditions hold for all $(i,j)$:

(1) $H_{ij} = 0 \Rightarrow H[i,j]$ is not sign-nonsingular
(2) $H_{ij} \neq 0 \Rightarrow H[i,j]$ is sign-nonsingular

**Proof:** The proposition follows directly from Lemma 2.3-1 and the fact that a non-zero entry in a sign-nonsingular pattern cannot correspond to a minor matrix whose determinant changes sign. ∎

We now show that for a sign-nonsingular pattern $H$ that allows symplectic pairs, the properties of maximality and irreducibility are equivalent. In particular, this implies that a strong noneven digraph $D$ that allows symplectic pairs must be maximal.

**Proposition 6.2:** Suppose $H$ is a negative-diagonal sign-nonsingular pattern that allows symplectic pairs. Then $H$ is maximal if and only if $H$ is irreducible.

**Proof:** First suppose $H$ is maximal. Then by the definition of maximality, $H_{ij} = 0 \Rightarrow H[i,j]$ is not combinatorially singular. Thus, given any pair of distinct indices $(i,j)$ such that $H_{ij} = 0$, there is at least one simple path from $j$ to $i$ in $D(H)$. Also, by Proposition 6.1, $H_{ij} \neq 0 \Rightarrow H[i,j]$ is sign-nonsingular, again guaranteeing at least one simple path from $j$ to $i$ in $D(H)$. Thus $D(H)$ is strongly connected and $H$ is irreducible.

Now suppose $H$ is irreducible. By Proposition 6.1, $H_{ij} = 0 \Rightarrow H[i,j]$ is not sign-nonsingular. Since $H$ is irreducible, $D(H)$ is strongly connected, so that $H$ contains no combinatorially singular minor matrix. It follows that if $H_{ij} = 0$, $H[i,j]$'s determinant changes sign, so that all zeros in $H$ are necessary zeros. $H$ is therefore maximal. ∎

We now characterize the strong, dense, noneven digraphs that allow symplectic pairs.

**Theorem 6.3:** Suppose $D$ is a strong, dense, noneven digraph. Then $D$ allows symplectic pairs if and only if $D = W_4$, $D = C_4^*$, or $D$ is an extended caterpillar.

**Proof:** As discussed in [13], $W_4$ allows symplectic pairs. Since $W_4$ and $C_4^*$ are sign-equivalent, and since the symplectic pair property is preserved under sign-equivalence, $C_4^*$ allows symplectic pairs as well. If $D$ is neither $W_4$ nor $C_4^*$, then by Theorems 5.1 and 5.2, $D$ is a subdigraph of an extended caterpillar. First suppose $D$ is an extended caterpillar. Then, since all extended caterpillars share the same indegree and outdegree sequences (both of which are $\{1, 2, 3, \ldots, n-1, n-1\}$), $D$ is sign-equivalent to the special caterpillar $C$, all of whose vertices are backbone vertices. As shown in [14], $C$ allows symplectic pairs, and so $D$ allows symplectic pairs as well. Now suppose $D$ is a proper subdigraph of an extended caterpillar. In this case, $D$ is not maximal, since $D$ was obtained by removing arcs from a noneven digraph. Thus $D$ cannot allow symplectic pairs by Proposition 6.2. ∎

# 7  Summary and Conclusions

Koh [15] showed that there are infinitely many noneven digraphs of minimum indegree and outdegree two. Thomassen showed [16] that there are infinitely many 2-connected noneven digraphs, and Lim [11] constructed an infinite family of sign-nonsingular patterns (associated with wheels) whose digraphs are 2-connected, noneven and maximal planar. Thomassen also showed [16] that no 3-connected digraph is noneven.

We have shown in this paper that among the dense 2-connected digraphs, the only noneven digraphs are of size four, namely $W_4$ and $C_4^*$, which are sign-equivalent. As shown in [1], dense digraphs of connectivity one must be subdigraphs of extended caterpillars. The only dense noneven digraphs that allow symplectic pairs are $W_4$, $C_4^*$, and extended caterpillars. (For detailed characterizations of the symplectic pairs allowed by the sign patterns corresponding to these digraphs, see [8] and [17].) More generally, no strong non-maximal noneven digraph can allow symplectic pairs.

Among the 2-connected symmetric digraphs, the only noneven digraphs cite1 are the digraphs $G^*$ where $G$ is a subdigraph of a $C_4$-cockade. The only such digraph that allows symplectic pairs is $C_4^*$ itself. The only symmetric digraph of connectivity one that allows symplectic pairs is $C_2^*$.